\documentstyle[12pt,amssymb,amsfonts,epsbox]{amsart}
\theoremstyle{plain}
\newtheorem{thm}{Theorem}[section]
\newtheorem{lem}[thm]{Lemma}

\newtheorem{prop}[thm]{Proposition}
\newcommand{\Diff}{\mathrm{Diff}}
\newcommand{\inter}{\mathrm{int}}
\newcommand{\fix}{\mathrm{fix}}
\newcommand{\Out}{\mathrm{Out}}
\theoremstyle{definition}

\theoremstyle{remark}
\newtheorem{rem}[thm]{Remark}

\begin{document}
\baselineskip=20pt
\title{Some homological invariants of mapping class group of 
a 3-dimensional handlebody} 
\author{Susumu Hirose}
%
\address{Department of Mathematics, 
Faculty of Science and Engineering, 
Saga University, 
Saga, 840 Japan}
\email{hirose@@ms.saga-u.ac.jp}
\subjclass{57N10, 57N05, 20F38}
\keywords{virtual cohomological dimension, Euler number, 
3-dimensional handlebody, mapping class group}
\date{February 16, 2001}
\maketitle
\begin{abstract}
We show that, 
if $g \geq 2$, the virtual cohomological dimension of the mapping class 
group of a 3-dimensional handlebody of genus $g$ is equal to $4g-5$ and 
the Euler number of it is equal to 0.  
\end{abstract}
\section{Introduction}
A genus $g$ {\it handlebody\/} $H_g$ is an oriented 3-manifold which is 
constructed from 3-ball with attaching $g$ 1-handles. 
The {\it mapping class group\/} ${\cal H}_{g}$ of $H_g$ is defined as the 
group of isotopy classes of orientation-preserving diffeomorphisms of $H_g$. 
This group ${\cal H}_{g}$ is a subgroup of the mapping class group 
${\cal M}_{g}$ of a surface $\partial H_g$, that is, 
${\cal M}_{g}$ $=$ $\pi_0( \Diff^{+} (\partial H_g))$, 
where $\Diff^{+}(\partial H_g)$ is the group of orientation preserving 
diffeomorphisms of $\partial H_g$. 
From here to the end of this paper, we assume $g \geq 2$. 
\par
The {\it cohomological dimension\/} of a group $G$, $cd(G)$, is defined to be 
the largest number $n$ for which there exist a $G$-module $M$ with 
$H^n(G,M)$ nonzero. 
We remark that if $G_1 \subset G_2$, then 
$cd(G_1)$ $\leq$ $ cd(G_2)$. 
When $G$ has torsion, $cd(G)$ is infinite. 
However, if $G$ has finite index torsion free subgroups (we call 
$G$ {\it virtually torsion free\/}), 
we define the {\it virtual cohomological dimension\/} of $G$, 
$vcd(G)$, to be the cohomological dimension of 
finite index torsion free subgroup $\hat{G}$. 
A theorem of Serre \cite{Serre} states that this number is independent 
of the choice of $\hat{G}$. 
For $vcd$ of ${\cal M}_{g}$ and 
${\cal H}_{g}$, Harer \cite{Harer} showed  $vcd({\cal M}_{g}) = 4g-5$, 
and  McCullough \cite{McCullough} showed 
$vcd({\cal H}_{2})$ $= 3$, and, if $g \geq 3$,  
$3g-2 \leq vcd({\cal H}_{g}) \leq 4g-5$. In this paper, we will show 
the following result. 
\begin{thm}\label{thm:vcd}
If $g \geq 2$, 
the virtual cohomological dimension of ${\cal H}_{g}$ is equal to 
$4g-5$. 
\end{thm}
Hatcher observed that this result can be shown by investigation of 
the action of ${\cal H}_{g}$ on the disc complex defined by 
McCullough \cite{McCullough}. 
In this paper, we show this result with the explicit description of 
subgroup of ${\cal H}_{g}$ that achieves the $vcd({\cal H}_{g})$.

We give some remark on the relationship between 
${\cal H}_{g}$ and the outer automorphism group of free group of rank $g$. 
We denote $F_g$ the free group of rank $g$ and 
$\Out(F_g)$ the outer automorphism group of it. 
There is a natural homomorphism from 
${\cal H}_{g}$ to $\Out(F_g)$ defined by the action of diffeomorphisms 
on the fundamental group of $H_g$. 
This homomorphism is a surjection \cite{Griffiths}. 
Culler and Vogtmann \cite{C-V} showed that $vcd(\Out(F_g))$ 
$=$ $2g-3$. 
This fact indicate that the kernel of the above surjection 
is, in some sense, big. 
In fact, McCullough \cite{McCullough2} showed that the kernel of 
the above surjection is not finitely generated.

We review the Euler characteristics of groups (see \cite{Brown}). 
For a group $G$ of finite homological type and torsion-free, 
we define the {\it Euler characteristic\/} $\chi(G)$ by 
$$
\chi(G) = \sum_i (-1)^i \text{rk}_{\Bbb{Z}} (H_i(G)). 
$$
For a group $G$ of finite homological type which may have torsion, 
we choose a torsion free subgroup $\hat{G}$ of finite index, and define 
$\chi(G)$ by 
$$
\chi(G) = \frac{\chi(G)}{(G:\hat{G})}, 
$$
where $(G:\hat{G})$ is the index of $\hat{G}$ in $G$. 
Since, ${\cal H}_{g}$ is of type VFL \cite{McCullough}, 
we can define $\chi({\cal H}_{g})$. 
We will show the following result. 
\begin{thm}\label{thm:Euler}
$\chi({\cal H}_{g}) = 0$. 
\end{thm}
\section{Proof of Theorem \ref{thm:vcd}.}
\begin{figure}
\begin{center}
\psbox[height=5cm]{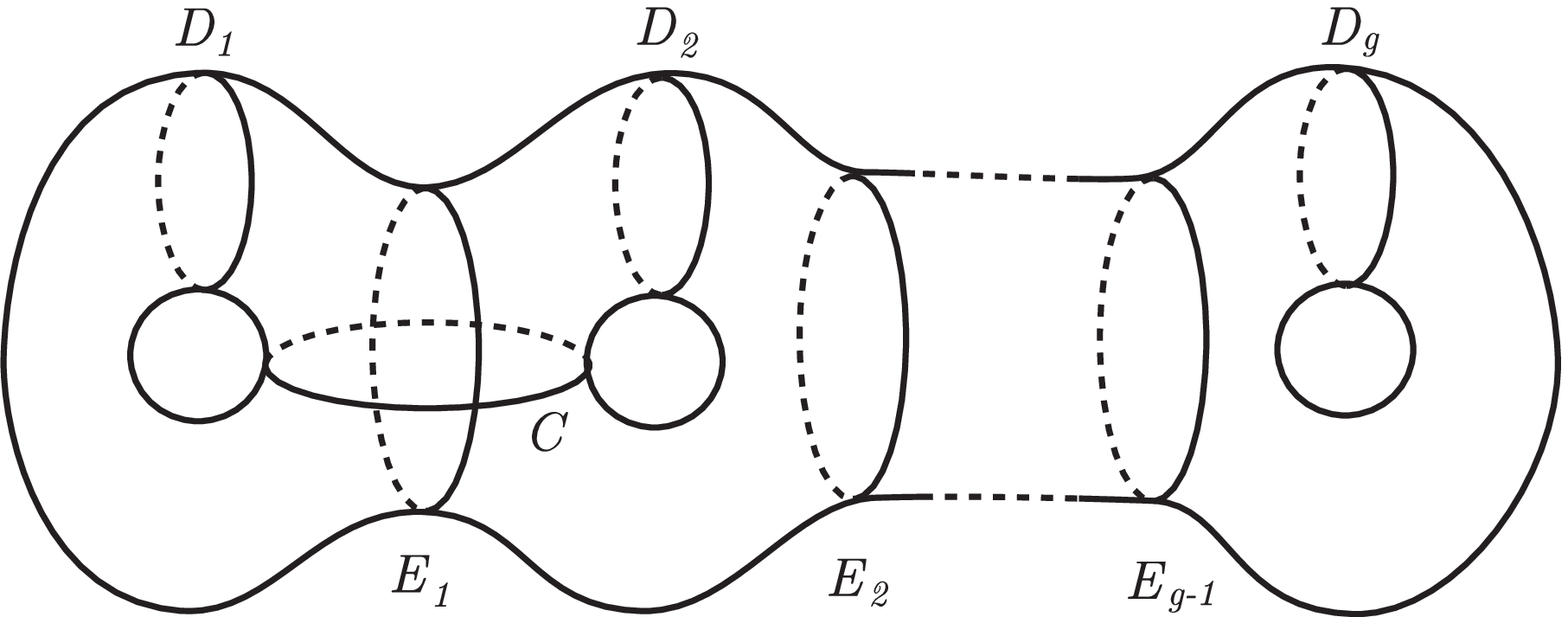}
\caption{}
\label{fig:discs}
\end{center}
\end{figure}

In general, for an oriented $C^{\infty}$-manifold $A$ and its subset 
$B$, we denote $\Diff^{+} (A)$ the group of 
all orientation preserving 
diffeomorphisms of $A$, 
denote $\Diff^+ (A,\ \fix \  B)$ the group of elements of 
$\Diff^{+} (A)$ whose restriction to $B$ is identity map, 
and denote $\Diff^+ (A,\  B)$ the group of elements of 
$\Diff^{+} (A)$ which fixes $B$ as a set. 
For a disc $D$ in $\partial H_g$, we define ${\cal H}_{g}^1$ $=$ 
$\pi_0 (\Diff^{+}(H_g, \fix \ D))$, and define ${\cal M}_{g}^1$ $=$ 
$\pi_0 (\Diff^{+}(\partial H_g, \fix \ D))$. 
For the center $p$ of the above disc $D$, we define 
${\cal H}_{g,1}$ $=$ $\pi_0 (\Diff^{+}(H_g, \fix \ \{p\}))$, 
and define ${\cal M}_{g,1}$ $=$ 
$\pi_0 (\Diff^{+}(\partial H_g, \fix \ \{p\}))$. 
Let $D_1, D_2, \ldots, D_g$ be the cocores of 1-handles which are used 
to construct $H_g$. These discs $D_1, D_2, \ldots, D_g$ are properly 
embedded discs in $H_g$. 
Let $E_1,\ldots, E_{g-1}$ and $C$ be properly embedded discs as are indicated 
in Figure \ref{fig:discs}.
\par
We introduce some elements of ${\cal H}_{g}$. 
For a disc $D$ properly embedded in $H_g$, let $N$ be a regular neighborhood 
of $D$ in $H_g$. 
We parametrize $N$ by $\phi:$ $[-1,1] \times D^2$ $\longrightarrow$ $N$ 
such that, $\phi(\{0\} \times D^2)$ $=$ $D$, and 
$\phi([-1,1] \times \partial D^2)$ is an annulus in $\partial H_g$. 
Let $\psi$ be the diffeomorphisms of $[-1,1]\times D^2$, 
defined by $\psi(t,r,\theta) = (t,r,\theta+(1-t) \pi)$, 
where $(r, \theta)$ is a polar coordinate of $D^2$. 
The map $\delta_D$ from  $H_g$ to itself, 
defined by $\delta_D(x)$ $=$ $\phi \circ \psi \circ \phi^{-1} (x)$ 
iff $x \in N$, $= x$ iff $x \not\in N$, 
is an orientation preserving diffeomorphism of $H_g$. 
We call this {\it disc twist} about $D$. 
The isotopy class of $\delta_D$ is an element of ${\cal H}_{g}$, 
and we call this {\it disc twist} about $D$ and denote this $d_D$. 
For an annulus $A$ properly embedded in $H_g$, let $N$ be a regular 
neighborhood of $A$ in $H_g$. 
We parametrize $N$ by $\phi:$ $[-1,1]\times S^1 \times [0,1]$ $\longrightarrow$ 
$N$, such that $\phi|_{\{0\} \times S^1 \times [0,1]}$ is a parametrization 
of $A$, $\phi([-1,1] \times S^1 \times \{0\})$ and 
$\phi([-1,1] \times S^1 \times \{1\})$ are disjoint annuli in 
$\partial H_g$. 
Let $\psi$ be the diffeomorphism on $[-1,1] \times S^1 \times [0,1]$ defined 
by $\psi(t, \theta, s) = (t, \theta + (1-t) \pi, s)$, where $\theta$ is a 
polar coordinate of $S^1$. 
The map $\alpha_A$ of $H_g$, defined by 
$\alpha_A(x)$ $=$ $\phi \circ \psi \circ \phi^{-1} (x)$ iff $x \in N$, 
$=$ $x$ iff $x \not\in N$, is an orientation-preserving diffeomorphism 
on $H_g$. 
We call this {\it annulus twist} about $A$.
The isotopy class of $\alpha_A$ is an element of ${\cal H}_{g}$, 
and we call this {\it annulus twist} about $A$ and denote this $a_A$. 
\par
We introduce some expressions used in this paper. 
Let $N$ be a regular neighborhood of $\partial H_g$ in $H_g$, and 
$A$ be an annulus in $\partial H_g$. 
We parametrize $N$ as 
$\phi:$ $[0,1]\times \partial H_g \longrightarrow N$, such that 
$\phi(\{0\} \times \partial H_g) = \partial H_g$ and 
$\phi|_{\{0\} \times \partial H_g}$ is an identity map. 
The set $A'$ $=$ $\phi(\partial A \times [0,1] \cup A \times \{1\})$ 
is an annulus properly embedded in $H_g$. 
The sentence "we {\it push\/} $A$ {\it into\/} $H_g$" means that 
we obtain $A'$ from $A$. 
For a disc $D$ in $\partial H_g$, the meaning of the sentence 
"we {\it push\/} $D$ {\it into\/} $H_g$" is given in the same manner as above.
\par
G.~Mess \cite{Mess} discovered some subgroups $B_g$, $B_g^1$ of the mapping 
class groups ${\cal M}_{g}$, ${\cal M}_{g}^1$ respectively. 
We call $B_g$ and $B_g^1$ {\it Mess subgroups\/}. 
We review the definition of Mess subgroups (this definition is quoted from 
\S 6.3 of \cite{Ivanov}). Mess subgroups are defined in a recursive manner. 
\par
{\it Step $0$\/}: Let $B_2$ be the subgroup of ${\cal M}_{2}$ generated by 
Dehn twist about any three pairwise disjoint pairwise nonisotopic 
simple closed curves $C_0, C_1, C_2$ in $\partial H_2$. 
\par
{\it Step $1_g$\/}: We assume that $B_g$ ($g \geq 2$) is already defined. 
There is a surjection from $\Diff^{+}(\partial H_g, \fix\ D)$ to 
$\Diff^{+}(\partial H_g)$ defined by forgetting the disc $D$, and this 
surjection induce an surjection $f:$ ${\cal M}_{g}^1 \longrightarrow {\cal M}_{g}$. 
Let $B_g^1$ be the preimage of $B_g$ under $f$. 
\par
{\it Step $2_g$\/}: By the restriction of diffeomorphisms, we can 
define a homomorphism $\rho:$ $\Diff^{+}(\partial H_g, \fix\ D)$ 
$\longrightarrow$ $\Diff^{+}(\partial H_g - \inter D, \fix \ \partial D)$. 
We consider some embedding $\partial H_g - \inter D$ into 
$\partial H_{g+1}$ and identify $\partial H_g - \inter D$ with its image. 
The extension of diffeomorphisms of $\partial H_g - \inter D$ fixed on 
$\partial D$ by the identity across the complement of 
$\partial H_g - \inter D$ in $\partial H_{g+1}$ define a homomorphism 
$\iota:$ $\Diff^{+} (\partial H_g - \inter D, \fix \partial D)$ 
$\longrightarrow$ $\Diff^{+} (\partial H_{g+1})$. 
The homomorphism $\iota \circ \rho$ induce a homomorphism 
$i:$ ${\cal M}_{g}^1 \longrightarrow {\cal M}_{g+1}$. 
In the complement of $\partial H_g - \inter D$ in $\partial H_{g+1}$, 
we choose some non-trivial simple close curve $C$ and consider the Dehn twist 
$t \in {\cal M}_{g+1}$ about this curve. 
Let $T$ be the infinite cyclic group generated by $t$. 
We define $B_{g+1}$ as the group generated by $i(B_g^1)$ and $T$. 
\par
Mess showed the following theorem \cite{Mess} (see also Corollary 6.3B 
of \cite{Ivanov}). 
\begin{thm}\label{thm:Mess}
The cohomological dimension of $B_g$ is equal to $4g-5$. 
\qed
\end{thm}
We will show the following lemma. 
\begin{lem}\label{lem:realisation}
$B_g$ is a subgroup of ${\cal H}_{g}$. 
\end{lem}
\begin{rem} 
The above fact is remarked by Mess \cite[p.4]{Mess}. 
\end{rem}
The definition of $B_g$ involves some choices. This lemma means 
that, with some good choices, $B_g$ is realised as a subgroup of ${\cal H}_{g}$. 
\par
\noindent
{\bf Proof. }
Along the steps of the definition of $B_g$, we will check that 
$B_2$, $B_g^1$, $B_{g+1}$ can be constructed as subgroups of 
${\cal H}_{2}$, ${\cal H}_{g}^1$, ${\cal H}_{g+1}$ respectively. 
In each steps, we use the same notations as used in definitions of 
$B_g$ and $B_g^1$. \par
{\it Step $0$\/}: We choose $C_0 = \partial D_1$, $C_1=\partial C$, 
$C_2 = \partial D_2$, then $B_2 \subset {\cal H}_{2}$. 
\par
{\it Step $1_g$\/}: We assume that $B_g$ $\subset$ ${\cal H}_{g}$. 
Let $g_1, \ldots, g_n$ be the generators of $B_g$. 
For each $g_i$, we can choose an element $\tilde{g}_i$ of ${\cal H}_{g}^1$ 
such that $f(\tilde{g}_i) = g_i$. 
By the definition, $B_g^1$ is generated by the kernel of $f$ and 
$\tilde{g}_1, \ldots, \tilde{g}_n$. 
In order to obtain generators for the kernel of $f$, we consider the following 
two short exact sequences. 
\begin{equation}
\tag{S1}
1 \longrightarrow {\Bbb Z}
\longrightarrow {\cal M}_{g}^1  
\overset{\alpha}{\longrightarrow} {\cal M}_{g,1} 
\longrightarrow 1. 
\end{equation}
\begin{equation}
\tag{S2}
1 \longrightarrow \pi_1 (\partial{H}_g, \fix \{p\}) 
\overset{\beta}{\longrightarrow} {\cal M}_{g}^1 
\overset{\gamma}{\longrightarrow} {\cal M}_{g} 
\longrightarrow 1. 
\end{equation}
The group ${\Bbb Z}$ in (S1) is an infinite cyclic group generated 
by the Dehn twist $d$ about $\partial D$. 
The homomorphism $\alpha$ is induced from the homomorphism from 
$\Diff^{+} (\partial H_g, \fix \ D)$ to 
$\Diff^{+} (\partial H_g, \fix \{ p \})$ defined by crushing $D$ into 
a point $p$. 
The sequence (S2) is introduced by Birman \cite{Birman}. 
The homomorphism $\gamma$ is induced from the homomorphism from 
$\Diff^{+} (\partial H_g, \fix \{ p \})$ to 
$\Diff^{+} (\partial H_g)$ defined by forgetting the point $p$. 
The group $\pi_0(\partial H_g, p)$ is generated by simple loops in 
$\partial H_g$, whose base points are $p$. 
Let $l_1, \ldots, l_{2g}$ be simple loops in $\partial H_g$, 
homotopy classes of which generates $\pi_0(\partial H_g, p)$. 
For each $l_i$, let $L_i$ be an annulus in $\partial H_g$, 
which is a regular neighborhood of $l_i$ such that $L_i \supset D \ni p$. 
$\partial L_i$ is two simple closed curves 
$l_i^1$, $l_i^2$ in $\partial H_g$. 
The homomorphism $\beta$ is defined so that 
it maps a homotopy class of $l_i$ (denote $[l_i]$ for short ) 
to a homotopy class of 
$\lambda_i$ $=$ $(+ \text{Dehn twist about } l_i^1) 
\times (- \text{Dehn twist about } l_i^2)$. 
This homeomorphism $\lambda_i$ is also an element of 
$\Diff^+ (\partial H_g, \fix \ D)$, and 
$\alpha ( $ an element of ${\cal M}_{g}^1$ represented by $\lambda_i)$ 
$=$ $\beta([l_i])$. 
Let $\tilde{l}_i$ be an element of ${\cal M}_{g}^1$ represented by 
$\lambda_i$. 
Since the kernel of $f$ is equal to $\alpha^{-1}($ the kernel of 
$\gamma)$ $=$ $\alpha^{-1}($ the image of $\beta)$, 
the kernel of $f$ is generated by $d$ and 
$\tilde{l}_1, \ldots, \tilde{l}_{2g}$. 
Let $D'$ be a disc in $H_g$ made 
by pushing $D$ into $H_g$, and $\delta_{D'}$ be the disc twist about $D'$. 
Let $L'_i$ be an annulus made by pushing $L_i$ into $H_g$, 
and $\alpha_{L'_i}$ be the annulus twist about $L'_i$. 
The diffeomorphisms $\delta_{D'}$ and $\alpha_{L'_i}$ are elements of 
$\Diff^+ (H_g, \fix \ D)$, and restrictions of them 
to $\partial H_g$ represent $d$ and $\tilde{l}_i$ respectively. 
This fact shows that the kernel of $f$ is included in ${\cal H}_{g}^1$. 
Hence, $B_g^1$ $\subset$ ${\cal H}_{g}^1$. 
\par
{\it Step $2_g$\/}: It is easy to see that 
$i(B_g^1) \subset {\cal H}_{g+1}$. 
We choose $C$ $=$ $\partial D_{g+1}$, then $t$ $\in$ ${\cal H}_{g+1}$. 
Therefore, $B_{g+1}$ $\subset$ ${\cal H}_{g+1}$. 
\qed

\par
Along the line of the proof of Theorem 6.4.A in \cite{Ivanov}, 
we will prove Theorem \ref{thm:vcd}. \par
\noindent
{\bf Proof of Theorem \ref{thm:vcd}.} \newline
There is a natural homomorphism ${\cal M}_{g}$ $\to$ 
$Aut(\mathrm{H}_1(\partial H_g,{\Bbb Z}/3{\Bbb Z}))$ defined by 
the action of diffeomorphisms on homology. 
Let $\Gamma$ be the kernel of this homomorphism. 
By Ivanov \cite[Corollary 1.5]{Ivanov-book}, $\Gamma$ is 
torsion free. 
Therefore, $\Gamma$, ${\cal H}_{g} \cap \Gamma$, and $B_g \cap \Gamma$ are 
finite index torsion free subgroups of 
${\cal M}_{g}$, ${\cal H}_{g}$, and $B_g$ respectively. 
By the definition of virtual cohomological dimension, 
$vcd({\cal M}_{g})$ $=$ $cd(\Gamma)$, 
$vcd({\cal H}_{g})$ $=$ $cd({\cal H}_{g} \cap \Gamma)$ and 
$vcd(B_g)$ $=$ $cd(B_g \cap \Gamma)$. 
By Harer \cite[Theorem 4.1]{Harer}, $vcd({\cal M}_{g})$ $=$ $4g-5$, hence, 
$cd(\Gamma)$ $=$ $4g-g$. 
By Theorem \ref{thm:Mess}, $vcd(B_g)$ $=$ $cd(B_g)$ $=$ $4g-5$, 
hence, $cd(B_g \cap \Gamma)$ $=$ $4g-5$. 
By Lemma \ref{lem:realisation}, $B_g \cap \Gamma$ $\subset$ 
${\cal H}_{g} \cap \Gamma$ $\subset$ $\Gamma$, therefore, 
$cd(B_g \cap \Gamma)$ $\leq$ $cd({\cal H}_{g} \cap \Gamma)$ 
$\leq$ $cd(\Gamma)$. 
These fact show this theorem. 
\qed
\section{Proof of Theorem \ref{thm:Euler}.}
McCullough defined a {\it disc complex\/} L in \cite{McCullough} and used 
this to give an estimation for $vcd({\cal H}_{g})$. 
We review the definition of $L$. 
By a {\it disc\/} in $H_g$, we mean a properly embedded 2-disc in $H_g$. 
The disc is {\it essential\/} when $\partial D$ does not bound a 2-disc in 
$\partial H_g$. 
The disc complex $L$ of $H_g$ is the simplicial complex whose vertices are 
the isotopy classes of essential discs in $H_g$, and whose simplices are 
defined by the rule that a collection of $n+1$ distinct vertices spans a 
$n$-simplex if and only if it admit a collection of representatives 
which are pairwise disjoint. 
McCullough showed the following Theorem. 
\begin{thm}\label{thm:disc-comp}\cite[Theorem 5.3]{McCullough}
The disc complex $L$ of $H_g$ is contractible. 
\qed
\end{thm}
We use the following two Propositions about Euler characteristics of groups. 
\begin{prop}\label{prop:exact}\cite[Proposition \S IX 7.3(d)]{Brown}
Let $1 \to G' \to G \to G'' \to 1$ be a short exact sequence of groups with 
$G'$ and $G''$ of finite homology type. 
If $G$ is virtually torsion free, then $G$ is of finite homological type and 
$\chi(G)$ $=$ $\chi(G') \chi(G'')$. 
\qed
\end{prop}
\begin{prop}\label{prop:equiv}\cite[Proposition \S IX 7.3(e')]{Brown}
Let $X$ be a contractible simplicial complex on which $G$ act simplicially. 
For each simplex $\sigma$ of $X$, let $G_{\sigma}$ $=$ 
$\{ g \in G | g \sigma = \sigma \}$. 
If $X$ has only finitely many cells mod $G$, and, for each simplex 
$\sigma$ of $X$, $G_{\sigma}$ is of finite homological type, 
then 
$$
\chi(G) = \sum_{\sigma \in \cal{E}} (-1)^{\dim \sigma} \chi(G_{\sigma}), 
$$
where $\cal{E}$ is a set of representative for the cells of $X$ mod $G$. 
\qed
\end{prop}
For each simplex $\sigma$ $=$ $<D_0, \ldots, D_n>$ of $L$, 
$G_{\sigma}$ $=$ $\pi_0 \Diff^+ (H_g,\ D_0\cup \cdots \cup D_n)$. 
For the same simplex $\gamma$, let $\Gamma_{\sigma}$ be the graph defined 
as follows. 
The vertices of $\Gamma_{\sigma}$ correspond to the components of 
$H_g - D_0 \cup \cdots \cup D_n$. 
Each edge corresponds to one of $D_0, \ldots, D_n$ and connects the vertices 
corresponding to the components attached along this disc. 
There is a natural homomorphism $\delta$ from $G_{\sigma}$ to the group of 
automorphisms of $\Gamma_{\sigma}$. 
Let $A_{\sigma}$ be the image of $\delta$. 
Since the group of automorphisms of $\Gamma_{\sigma}$ is a finite group, 
$A_{\sigma}$ is a finite group. 
There are short exact sequences. 
\begin{equation}
\tag{S3}
1 \longrightarrow \pi_0 \Diff^+ (H_g,\ \fix \ D_0 \cup \cdots \cup D_n) 
\longrightarrow G_{\sigma} 
\overset{\delta}{\longrightarrow} A_{\sigma} 
\longrightarrow 1. 
\end{equation}
\begin{align*}
\tag{S4}
1 \longrightarrow & {\Bbb Z}^{n+1}
\longrightarrow \pi_0 \Diff^+ (H_g,\ \fix \ D_0 \cup \cdots \cup D_n) \\
&\overset{\epsilon}{\longrightarrow} 
\pi_0 \Diff^+ (H_g/D_0 \cup \cdots \cup D_n,\ \fix 
\ D_0/D_0 \cup \cdots \cup D_n/D_n) 
\longrightarrow 1. 
\end{align*}
The homomorphism $\epsilon$ is induced by crashing each $D_i$ into one point. 
The group ${\Bbb Z}^{n+1}$ is generated by disc twist about $D_0$, $\cdots$, 
$D_n$, and , as is well-known, $\chi({\Bbb Z}^{n+1}) = \chi ((S^1)^{n+1}) = 0$. 
By applying Proposition \ref{prop:exact} to (S4), 
we see $\chi(\pi_0 \Diff^+ (H_g,\ \fix \ D_0 \cup \cdots \cup D_n)) = 0$. 
Therefore, by (S3), we obtain $\chi(G_{\sigma})=0$. 
Theorem \ref{thm:Euler} follows from the above observation, 
Theorem \ref{thm:disc-comp} and Proposition \ref{prop:equiv}. 

\begin{center}
{\sc Acknowledgements}
\end{center}
The author would like to express his gratitude 
to Prof. T.~Akita and Prof. N.~Kawazumi and 
Prof. D.~McCullough 
for their helpful comments.

\end{document}